\documentclass[12pt]{article}
\usepackage{latexsym}
\usepackage{amsfonts}
\def\ZZ{\mathbb Z}
\def\RR{\mathbb R}

\def\CP{\mathbb C \mathbb P}

\def\bea{\begin{eqnarray*}}
\def\eea{\end{eqnarray*}}

\newtheorem{thm}{Theorem}

\newenvironment{proof}{\medskip \noindent
{\bf Proof.}}{\hfill \rule{.5em}{1em}
\\}

\begin{document}
\sloppy
\title{Scalar Curvature,  Covering Spaces, and Seiberg-Witten Theory}

\author{Claude LeBrun\thanks{Supported 
in part by  NSF grant DMS-0072591.}  
  }

\date{}

\maketitle

\begin{abstract} 
The {\em Yamabe invariant} ${\mathcal Y}(M)$
of a 
smooth compact manifold is roughly  the supremum of 
the scalar curvatures of unit-volume constant-scalar-curvature
Riemannian metrics $g$ on $M$. (To be 
precise, one only considers  those constant-scalar-curvature 
metrics which are {\sl Yamabe minimizers}, but this technicality  
does not, e.g.   affect the {\em sign} of the answer.)
In this article, it is shown that many   
 $4$-manifolds $M$ with ${\mathcal Y}(M) < 0$ have 
 have finite covering spaces $\tilde{M}$ 
 with ${\mathcal Y}(\tilde{M}) > 0$. 
 \end{abstract}

Two decades ago,  Lionel B\'erard Bergery \cite{beber} pointed out that 
there are high-dimensional smooth compact manifolds $M$  which   do not 
admit metrics of positive scalar curvature, but which nevertheless
have finite coverings that {\em do} admit such metrics. 
For example, let $\Sigma$ be an exotic $9$-sphere which
does not bound a spin manifold,  and consider the connected sum
$M=(S^2 \times {\mathbb R \mathbb P}^7) \# \Sigma .$
On one hand,   $M$  is a spin manifold
with  non-zero
Hitchin invariant $\hat{\mathfrak a}(M)\in \ZZ_2$, so \cite{hitharm}
there are harmonic spinors on $M$ for every choice of metric;
 the  Lichnerowicz  Weitzenb\"ock formula for the Dirac operator therefore
tells us that no metric on $M$ can have 
 positive scalar curvature. On the other hand, 
the universal cover $\tilde{M}=(S^2\times S^7)\# 2\Sigma $ of $M$ is diffeomorphic to 
$S^2\times S^7$, on which the obvious product metric  certainly has   positive
scalar curvature.

As  will be shown here, the   same phenomenon  also occurs in 
dimension four. Indeed, far more is  true. 
In the process of passing from a 
$4$-manifold to a finite 
  cover, it is even possible to change the 
 sign of the {\em Yamabe invariant}.

 The Yamabe invariant 
 is a diffeomorphism invariant that  
 historically arose from  an  attempt to construct Einstein metrics 
 (metrics of constant Ricci curvature) 
on  smooth compact 
manifolds. A standard computation  \cite{bes} shows  that the 
Einstein metrics on any given smooth compact manifold $M$ of dimension $n > 2$    are 
exactly  the critical points of 
the {\em normalized total scalar curvature}
$${\mathfrak   S} (g)=V_g^{(2-n)/{n}} \int_{M}s_{g}d\mu_{g} ,$$
considered as a functional on the space of all Riemannian metrics $g$ on
$M$; 
here   $s$, $d\mu$, and $V$ respectively  denote the scalar curvature,  
   volume measure, 
and  total 
volume  of the relevant metric. 
However, one cannot possibly hope to 
find a critical point of  ${\mathfrak   S}$  by either  maximizing or minimizing  it, as it  is 
 bounded neither above nor  below. Nevertheless, as was first 
pointed out  by Hidehiko Yamabe \cite{yam},
 the restriction of 
${\mathfrak  S}$ to any {\em conformal class}
$$\gamma = [g] = \{ ug ~|~ u: M \stackrel{C^{\infty}}{\longrightarrow}
{\RR}^{+}\}$$ 
of metrics {\em  is} always 
bounded below. The trail blazed by Yamabe eventually led 
 \cite{trud,aubyam,rick} 
 to  a proof   of the fact that, for each conformal class $\gamma$, this
 infimum is actually achieved,  by a constant-scalar-curvature metric
known as a {\em Yamabe minimizer}. 
Yamabe's ultimate goal  was   to  construct   
 Einstein metrics by  maximizing the restriction of $\mathfrak S$ to the 
 set of these Yamabe minimizers. 
This last idea turns out to be  unworkable in practice, 
 but it nonetheless  gives rise to a beautiful,  
real-valued 
diffeomorphism invariant  \cite{okob,sch,lno}
$${\mathcal Y}(M)=
\sup_{\gamma}\inf_{g\in\gamma}
{\mathfrak  S}(g),$$
 called the {\em Yamabe invariant}  (or {\em sigma constant}) of  $M$. 
It is not hard to show that 
  ${\mathcal Y}(M) > 0$
if and only if 
 $M$ admits a metric of positive scalar curvature; thus 
the problem of computing the Yamabe invariant 
may be thought of as a quantitative refinement 
of the question of whether a given manifold admits 
positive-scalar-curvature metrics. 
On the other hand, 
if $M$ does {\em not} admit metrics of positive scalar curvature, 
the number 
${\mathcal Y}(M) \leq  0$ is just  the   supremum of 
the scalar curvatures of all unit-volume 
constant-scalar-curvature metrics on $M$. Dimension
$4$ turns out to be remarkably special  so far as this invariant is concerned. Indeed, 
Seiberg-Witten theory allows one to show  \cite{lno} that there are many simply connected 
$4$-manifolds with ${\mathcal Y}(M) < 0$. By contrast, however,  Petean \cite{jp3} has  shown
 that  every
simply connected compact manifold $M^n$ of dimension $n\geq 5$ 
has ${\mathcal Y} (M) \geq 0$. 

The  main construction used in this paper 
 primarily depends  on properties of the oriented $4$-manifold
$$N= (S^2\times S^2) /\ZZ_2 , $$
where the $\ZZ_2$ action is generated by 
the double antipodal map 
$$(\vec{x},\vec{y})\mapsto (-\vec{x},-\vec{y}).$$
Let   $X$ be any   non-spin compact complex surface of 
general type which can be expressed  as a  
{\em complete intersection} of complex hypersurfaces   in some complex 
projective space; for example, one could 
take $X$ to be the Fermat hypersurface
$$\Big\{ [v:w:x:y]\in {\CP}_3 ~|~    v^m + w^m + x^m + y^m =0\Big\}$$ 
for any odd $m\geq 5$.

\begin{thm} \label{first} 
Let $X$ and $N$ be as above, and let $M=X\# N$. 
Then $M$ has negative Yamabe invariant, but its universal cover $\tilde{M}$
  has positive Yamabe invariant. 
\end{thm}
\begin{proof}
Since  $X$ is a complex algebraic surface with 
$b_+(X) > 1$, the Seiberg-Witten invariant of $X$ is 
well-defined and non-zero for the canonical spin$^c$ structure
determined by the complex structure. 
 On the other hand, 
$N$ satisfies $b_2(N)=b_1(N)=0$, and a  gluing  result  
of Kotschick-Morgan-Taubes \cite{kmt} thus implies
that the Seiberg-Witten invariant is non-zero for 
the associated  spin$^c$ structure on $M=X\# N$ with $c_1=c_1(X)$.
This tells us \cite{witten}  that $M$ does not admit any metrics of positive
scalar curvature, and that \cite{lno}, moreover,
$${\mathcal Y}(X\# N) \leq -4\pi \sqrt{2c_1^2(X)}.$$
In particular,  ${\mathcal Y}(X\# N) <0$.

On the other hand, the universal cover of $M$ is
$\tilde{M}= X\# X\# (S^2 \times S^2)$. 
But   
 Gompf \cite{gfsum},
  inspired by the earlier work of  Mandelbaum and Moishezon  \cite{mamoi}, 
 has used a handle-slide argument  to show  that  $X\#(S^2 \times S^2)$ {\em dissolves},
in the sense that 
$$ X\#(S^2 \times S^2)\stackrel{\rm \tiny diff}{\approx} k_1 \CP_2 \#\ell_1
\overline{\CP}_2, $$
where $k_1 > 2$ and $\ell_1 \neq 0$. 
Since $(S^2 \times S^2)\# \CP_2 \approx 2 \CP_2 \# \overline{\CP}_2,$
it  follows that 
$$
X\# 2 \CP_2 \# \overline{\CP}_2 \approx (k_1+1) \CP_2 \#\ell_1
\overline{\CP}_2, $$
and hence that 
$$\tilde{M} = 2X \#(S^2 \times S^2)\approx 
X \# k_1 \CP_2 \#\ell_1
\overline{\CP}_2
\approx 
k \CP_2 \#\ell
\overline{\CP}_2, $$
where $k=2b_+(X) +1$ and $\ell=2 b_-(X) +1$. 
Since a connected sum of positive-scalar-curvature manifolds 
admits metrics of positive scalar curvature \cite{gvln,syrger}, we thus conclude 
that ${\mathcal Y}(\tilde{M}) > 0$. \end{proof}

It is unclear whether an analogous change in the sign of the 
Yamabe invariant ever occurs in higher dimensions. At any rate, 
 this phenomenon certainly {\em does not} occur
in  B\'erard Bergery's examples. For example,  we obviously have 
${\mathcal Y}(S^2\times {\mathbb R
\mathbb P}^7) > 0$;  and we also know that ${\mathcal Y}(\Sigma^9)=0$ 
by Petean's theorem \cite{jp3}. The Petean-Yun surgery theorem \cite{petyun} therefore implies that
their connected sum has ${\mathcal Y} \geq 0$, too; and   since $(S^2\times {\mathbb R
\mathbb P}^7)\# \Sigma^9$ does not admit metrics of positive scalar curvature, this 
 shows  that 
 ${\mathcal Y}([S^2\times {\mathbb R
\mathbb P}^7]\# \Sigma)=0$. 
(Indeed, so far as  we seem to know at present, {\em every} compact $n$-manifold with $|\pi_1| < \infty$
and 
$n \geq 5$ could turn out to have non-negative Yamabe invariant; for an interesting
 partial result in this direction, 
see  \cite{botros}.) 

It is also perhaps worth mentioning that one can actually  
 compute  the exact value of the Yamabe invariant for any of 
the manifolds $M=X\# N$ considered in 
  Theorem \ref{first}. Indeed, 
as already noted, the Seiberg-Witten argument tells us that
${\mathcal Y}(M)\leq  -4\pi \sqrt{2c_1^2(X)}$. 
On the other hand, 
 ${\mathcal Y}(N ) > 0$  and ${\mathcal Y}(X) < 0$,  a general inequality due to 
Osamu Kobayashi  \cite{okob} tells us that 
$${\mathcal Y}(X\# N)\geq {\mathcal Y}(X). 
$$
However, because $X$ is a minimal complex surface of general type, 
its Yamabe invariant is given \cite{lno}  by 
 $ {\mathcal Y}(X)= -4\pi \sqrt{2c_1^2(X)}$. 
The above inequalities therefore  allow us to ascertain the   exact value 
$${\mathcal Y}(X\# N)  =-4\pi \sqrt{2c_1^2(X)}$$
of the  Yamabe invariant for any of the manifolds in question. 

By contrast, however,  
 exact calculations of the Yamabe invariant
are notoriously  difficult in the positive case, owing to the fact that 
in the positive  regime  a
 constant-scalar-curvature metric need not be a Yamabe minimizer. 
 However, we do know 
\cite{lcp2}
 that ${\mathcal Y}(\CP_2) = {\mathcal Y}(\overline{\CP}_2) =
12\pi \sqrt{2}$. Thus Kobayashi's inequality  \cite{okob} 
predicts that any connected sum  of $\CP_2$'s and 
$\overline{\CP}_2$'s satisfies 
$${\mathcal Y}(
k \CP_2 \#\ell
\overline{\CP}_2 ) \in [ {\mathcal Y}(\CP_2 ) , {\mathcal Y}(S^4) ] = [12\pi\sqrt{2} , 8 \pi \sqrt{6}] ,
$$
and we can thus at least  conclude that  the Yamabe invariant 
of the corresponding  universal cover $\tilde{M}$ is always somewhere in 
 this narrow range. 

Finally, let us observe that 
the examples in 
Theorem \ref{first} can  be greatly generalized,  provided one does not
insist on passing to the {\em universal} cover.
 
\begin{thm}
Let $Y$ be a  symplectic 4-manifold with $b_+ > 1$, 
 $|\pi_1| < \infty$, and 
  non-spin universal cover.
 Let $N=(S^2\times S^2)/\ZZ_2$,
as before. Then $M=Y\# 2N$ does not admit metrics of
positive scalar curvature, but nonetheless has finite coverings
$\tilde{M}\to M$ which {\em do} carry such  metrics. 
Moreover, if the symplectic minimal model  of $Y$ has $c_1^2  \neq 0$,
the Yamabe invariant reverses sign as one passes from $M$ to $\tilde{M}$. 
\end{thm}

\begin{proof}
By a celebrated theorem of Taubes \cite{taubes},  the canonical 
spin$^c$ structure of the symplectic manifold
 $Y$ has non-zero Seiberg-Witten invariant, and
since  
$N\# N$ has $b_1=b_2=0$,  the same gluing argument \cite{kmt} as before implies
that the Seiberg-Witten invariant is non-zero for 
a spin$^c$ structure on $M=Y\# 2N$ with $c_1=c_1(Y)$.
It thus follows that $M$ does not admit metrics of
positive scalar curvature.  But even more is true.  By another remarkable  result of Taubes
\cite{taubes3},  we can express $Y$ as an iterated  symplectic blow-up of 
a symplectic manifold $Y_0$, called the {\em symplectic minimal model of}
$Y$,  which contains no symplectic $(-1)$-spheres,
and  satisfies $c_1^2 (Y_0) \geq 0$. On the other hand, 
 the same argument used in  the proof  of \cite[Theorem 2]{lno} then shows that, 
 for every metric $g$ on $M$, there is a Seiberg-Witten basic class for which 
$(c_1^+)^2 \geq c_1^2 (Y_0)$. The estimate 
${\mathcal Y}(M) \leq -4\pi \sqrt{2 c_1^2 (Y_0 )}$ then follows immediately. 
This shows that  $M$ actually has 
negative Yamabe invariant whenever  $c_1^2(Y_0)\neq 0$.

Next, let $X$ denote the universal cover of $Y$, and 
observe that $Y\# 2N$ has an $n$-fold cover of the
form $X\# 2nN$, where $n=|\pi_1(Y)|$. 
Thus, unfolding copies of $N$ one by one, we  obtain  a sequence of 
covering spaces  $\tilde{M}_\ell\to M$ 
with  
$$\tilde{M}_\ell = 2^\ell X~\# ~\left(2^\ell -1\right) (S^2\times S^2)~ \# 
~ \left[2^{\ell+1}(n-1)+2\right] N$$
 for each $\ell \geq 1$.
On the other hand, since $X$ is simply connected and non-spin by 
assumption, a justly famous result of Wall \cite{wall}  asserts that 
 there is an integer $k_0$ such that 
$$X \# k_0  (S^2 \times S^2) \approx [k_0 + b_+(X) ]\CP_2
\# [k_0 + b_-(X)] \overline{\CP}_2. $$
Since 
$(S^2\times S^2) \# \CP_2   \approx 
2\CP_2 \# \overline{\CP}_2$, it then  follows that 
$$X \#( k_0+1)\CP_2\# k_0 \overline{\CP}_2 
 \approx [k_0+1 + b_+(X) ]\CP_2
\# [k_0 + b_-(X)] \overline{\CP}_2 . $$
Induction therefore gives us 
$$m X \#  k (S^2 \times S^2) \approx [k  + mb_+(X) ]\CP_2
\# [k + mb_-(X)] \overline{\CP}_2 $$
for any $k\geq k_0$ and $m\geq 1$. Thus, setting $\tilde{M} = \tilde{M}_\ell$ for any 
 $\ell > \log_2 k_0$, we have  constructed a 
  finite covering 
$\tilde{M} \to M$ with 
$$\tilde{M} \approx  p \CP_2
 \# q\overline{\CP}_2 \#r N,$$
and  such  an $\tilde{M}$ admits positive-scalar-curvature metrics because
it is a connected sum of manifolds with positive scalar curvature. 
\end{proof}

It is perhaps  worth noting that when the  given $Y$  is  not simply-connected, 
essentially the same argument would also work for $Y \# N$.
\bigskip 

\noindent {\bf Acknowledgment.} The author would like to warmly thank
 Ming Xu for drawing his attention to the problem.


\begin{thebibliography}{10}

\bibitem{aubyam}
{\sc T.~Aubin}, {\em \'{E}quations diff\'erentielles non lin\'eaires et
  probl\`eme de {Y}amabe concernant la courbure scalaire}, J. Math. Pures Appl.
  (9), 55 (1976), pp.~269--296.

\bibitem{beber}
{\sc L.~B\'erard~Bergery}, {\em Scalar curvature and isometry group}, in
  Spectra of Riemannian Manifolds, Tokyo, 1983, Kagai Publications, pp.~9--28.

\bibitem{bes}
{\sc A.~Besse}, {\em {E}instein Manifolds}, Springer-Verlag, 1987.

\bibitem{botros}
{\sc B.~Botvinnik and J.~Rosenberg}, {\em The {Y}amabe invariant for non-simply
  connected manifolds}.
\newblock e-print, math.DG/0104186.

\bibitem{gfsum}
{\sc R.~E. Gompf}, {\em On sums of algebraic surfaces}, Invent. Math., 94
  (1988), pp.~171--174.

\bibitem{gvln}
{\sc M.~Gromov and H.~B. Lawson}, {\em The classification of simply connected
  manifolds of positive scalar curvature}, Ann. Math., 111 (1980),
  pp.~423--434.

\bibitem{hitharm}
{\sc N.~Hitchin}, {\em Harmonic spinors}, Advanecs in Mathematics, 14 (1974),
  pp.~1--55.

\bibitem{okob}
{\sc O.~Kobayashi}, {\em Scalar curvature of a metric of unit volume}, Math.
  Ann., 279 (1987), pp.~253--265.

\bibitem{kmt}
{\sc D.~Kotschick, J.~W. Morgan, and C.~H. Taubes}, {\em Four-manifolds without
  symplectic structures but with nontrivial {S}eiberg-{W}itten invariants},
  Math. Res. Lett., 2 (1995), pp.~119--124.

\bibitem{lno}
{\sc C.~LeBrun}, {\em Four-manifolds without {E}instein metrics}, Math. Res.
  Lett., 3 (1996), pp.~133--147.

\bibitem{lcp2}
\leavevmode\vrule height 2pt depth -1.6pt width 23pt, {\em Yamabe constants and
  the perturbed {S}eiberg-{W}itten equations}, Comm. An. Geom., 5 (1997),
  pp.~535--553.

\bibitem{mamoi}
{\sc R.~Mandelbaum and B.~Moishezon}, {\em On the topology of simply connected
  algebraic surfaces}, Trans. Amer. Math. Soc., 260 (1980), pp.~195--222.

\bibitem{jp3}
{\sc J.~Petean}, {\em The {Y}amabe invariant of simply connected manifolds}, J.
  Reine Angew. Math., 523 (2000), pp.~225--231.

\bibitem{petyun}
{\sc J.~Petean and G.~Yun}, {\em Surgery and the {Y}amabe invariant}, Geom.
  Funct. Anal., 9 (1999), pp.~1189--1199.

\bibitem{rick}
{\sc R.~Schoen}, {\em Conformal deformation of a {R}iemannian metric to
  constant scalar curvature}, J. Differential Geom., 20 (1984), pp.~478--495.

\bibitem{sch}
\leavevmode\vrule height 2pt depth -1.6pt width 23pt, {\em Variational theory
  for the total scalar curvature functional for {R}iemannian metrics and
  related topics}, Lec. Notes Math., 1365 (1987), pp.~120--154.

\bibitem{syrger}
{\sc R.~Schoen and S.~T. Yau}, {\em On the structure of manifolds with positive
  scalar curvature}, Manuscripta Math., 28 (1979), pp.~159--183.

\bibitem{taubes}
{\sc C.~H. Taubes}, {\em The {S}eiberg-{W}itten invariants and symplectic
  forms}, Math. Res. Lett., 1 (1994), pp.~809--822.

\bibitem{taubes3}
\leavevmode\vrule height 2pt depth -1.6pt width 23pt, {\em The
  {S}eiberg-{W}itten and {G}romov invariants}, Math. Res. Lett., 2 (1995),
  pp.~221--238.

\bibitem{trud}
{\sc N.~Trudinger}, {\em Remarks concerning the conformal deformation of
  metrics to constant scalar curvature}, Ann. Scuola Norm. Sup. Pisa, 22
  (1968), pp.~265--274.

\bibitem{wall}
{\sc C.~T.~C. Wall}, {\em On simply connected 4-manifolds}, J. Lond. Math.
  Soc., 39 (1964), pp.~141--149.

\bibitem{witten}
{\sc E.~Witten}, {\em Monopoles and four-manifolds}, Math. Res. Lett., 1
  (1994), pp.~809--822.

\bibitem{yam}
{\sc H.~{Y}amabe}, {\em On the deformation of {R}iemannian structures on
  compact manifolds}, Osaka Math. J., 12 (1960), pp.~21--37.

\end{thebibliography}
  \end{document}